# DIXMIER TRACE AND THE FOCK SPACE

HÉLÈNE BOMMIER-HATO, MIROSLAV ENGLIŠ and EL-HASSAN YOUSSFI

ABSTRACT. We give criteria for products of Toeplitz and Hankel operators on the Fock (Segal-Bargmann) space to belong to the Dixmier class, and compute their Dixmier trace. At the same time, analogous results for the Weyl pseudodifferential operators are also obtained.

## 1. INTRODUCTION

Consider the Segal-Bargmann-Fock space $\mathcal{F}_\gamma = \mathcal{F}_\gamma(\mathbf{C}^n)$ of all entire functions on $\mathbf{C}^n$ square-integrable with respect to the scaled Gaussian

$$d\mu_\gamma(z) := e^{-\gamma|z|^2}\,dz, \qquad \gamma > 0,$$

where $dz$ stands for the Lebesgue measure on $\mathbf{C}^n$. The function

$$K_\gamma(z,w) := \Big(\frac{\gamma}{\pi}\Big)^n e^{\gamma\langle z,w\rangle}, \qquad z,w \in \mathbf{C}^n,$$

is the reproducing kernel of $\mathcal{F}_\gamma$:

$$f(z) = \int_{\mathbf{C}^n} f(w) K_\gamma(z,w)\,d\mu_\gamma(w), \qquad \forall f \in \mathcal{F}_\gamma, \forall z \in \mathbf{C}^n,$$

and at the same time serves as the integral kernel

$$P_\gamma f(z) = \int_{\mathbf{C}^n} f(w) K_\gamma(z,w)\,d\mu_\gamma(w), \qquad z \in \mathbf{C}^n,$$

of the orthogonal projection $P_\gamma : L^2(\mathbf{C}^n, d\mu_\gamma) \equiv L^2_\gamma \to \mathcal{F}_\gamma$. For a bounded measurable function $f$ on $\mathbf{C}^n$, the Toeplitz and the Hankel operator with symbol $f$ are defined by

$$\begin{aligned} &T_f : \mathcal{F}_\gamma \to \mathcal{F}_\gamma, && T_f g = P_\gamma(fg); \\ &H_f : \mathcal{F}_\gamma \to L^2_\gamma \ominus \mathcal{F}_\gamma, && H_f g = (I - P_\gamma)(fg), \end{aligned} \tag{1}$$

respectively. These operators have been studied extensively owing to their relevance in operator theory, quantum mechanics and time-frequency analysis [4] [5] [13] [7]

1991 *Mathematics Subject Classification.* Primary 47B35; Secondary 47G30, 47B10, 30H20.

*Key words and phrases.* Fock space, Weyl calculus, Toeplitz operator, Hankel operator, Dixmier trace.

Research supported by Barrande grant MEB021108, AV ČR research plan AV0Z10190503, GA ČR grant 201/09/0473, and by the French ANR DYNOP, Blanc07-198398.

[24]. It is easily seen that $\|T_f\| \le \|f\|_\infty$ and that $T_f$ is compact if $f(z) \to 0$ as $|z| \to +\infty$, and similarly for $H_f$; more accurate criteria for compactness and $p$-Schatten class membership are also available [10] [2], though as of this writing there seem to be known no if-and-only-if conditions for the latter if $p \neq 2$.

Toeplitz and Hankel operators have long been studied also in the context of $\mathcal{F}_\gamma$ replaced by the weighted Bergman spaces $A^2_\alpha = \{f \in L^2(\mathbf{B}^n, (1-|z|^2)^\alpha \, dz) : f \text{ is holomorphic on } \mathbf{B}^n\}$, $\alpha > -1$, on the unit ball $\mathbf{B}^n$ of $\mathbf{C}^n$. It has been known for some time that for $f$ holomorphic and $n > 1$, the Hankel operator $H_{\overline{f}}$ then belongs to the Schatten class $\mathcal{S}^p$ if and only if $f$ is in the diagonal Besov space $B^p(\mathbf{B}^n)$ and $p > 2n$, or $f$ is constant and $p \le 2n$; this phenomenon is called a *cutoff* at $p = 2n$ [1]. In dimension $n = 1$, the situation is slightly different in that the cutoff occurs not at $p = 2$ but at $p = 1$. Since it is immediate from (1) that for holomorphic functions $f$ and $g$,

$$[T_{\overline{f}}, T_g] = T_{\overline{f}g} - T_g T_{\overline{f}} = H^*_{\overline{g}} H_{\overline{f}}, \tag{2}$$

one can rephrase the results just mentioned also in terms of $\mathcal{S}^p$-membership of the commutators $[T_{\overline{f}}, T_g]$. In any case, it follows that there are no non-zero trace-class Hankel operators $H_{\overline{f}}$, $f$ holomorphic, on $\mathbf{B}^n$ if $n = 1$, and similarly the product $H^*_{\overline{f}_1} H_{\overline{f}_2} \dots H^*_{\overline{f}_{2n-1}} H_{\overline{f}_{2n}} = [T_{\overline{f}_2}, T_{f_1}] \dots [T_{\overline{f}_{2n}}, T_{f_{2n-1}}]$, with $f_1, \dots, f_{2n}$ holomorphic, is never trace-class if $n > 1$. In particular, there is no hope for $n > 1$ of having an analogue of the well-known formula for the unit disc $\mathbf{D} = \mathbf{B}^1$,

$$\operatorname{tr}[T_{\overline{f}}, T_f] = \frac{1}{\pi} \int_{\mathbf{D}} |f'(z)|^2 \, dz, \tag{3}$$

expressing the trace of the commutator $[T_{\overline{f}}, T_f]$, with $f$ holomorphic, as the square of the Dirichlet norm. A remarkable substitute for (3) on $\mathbf{B}^n$ is the result of Helton and Howe [16], who showed that for smooth functions $f_1, \dots, f_{2n}$ on the closed ball, the complete anti-symmetrization $[T_{f_1}, T_{f_2}, \dots, T_{f_{2n}}]$ of the $2n$ operators $T_{f_1}, \dots, T_{f_{2n}}$ is trace-class and

$$\operatorname{tr}[T_{f_1}, \dots, T_{f_{2n}}] = \int_{\mathbf{B}^n} df_1 \wedge df_2 \wedge \dots \wedge df_{2n}.$$

There is, however, a generalization of (3) to $\mathbf{B}^n$, $n > 1$, in a different direction — using the Dixmier trace. This is notable especially in view of the prominent applications of the latter in Connes' noncommutative geometry [8]. Namely, it was shown by Guo, Zhang and one of the present authors [11] that for $f_1, \dots, f_n$ and $g_1, \dots, g_n$ smooth on the closed ball, the product $T = [T_{f_1}, T_{g_1}] \dots [T_{f_n}, T_{g_n}] = (T_{f_1} T_{g_1} - T_{g_1} T_{f_1}) \dots (T_{f_n} T_{g_n} - T_{g_n} T_{f_n})$ belongs to the Dixmier ideal $\mathcal{S}^{\mathrm{Dixm}}$ and has Dixmier trace equal to

$$\operatorname{tr}_\omega T = \frac{1}{n!} \int_{\mathbf{S}^{2n-1}} \prod_{j=1}^n \{f_j, g_j\}_b \, d\sigma, \tag{4}$$

where $d\sigma$ is the normalized surface measure on the boundary $\mathbf{S}^{2n-1}$ of $\mathbf{B}^n$ and $\{f, g\}_b$ is a "boundary Poisson bracket" given by

$$\{f, g\}_b = \sum_{j=1}^n \Big( \frac{\partial f}{\partial \overline{z}_j} \frac{\partial g}{\partial z_j} - \frac{\partial f}{\partial z_j} \frac{\partial g}{\partial \overline{z}_j} \Big) - (\overline{R} f R g - R f \overline{R} g), \tag{5}$$

with

$$\overline{R} := \sum_{j=1}^{n} \overline{z}_j \frac{\partial}{\partial \overline{z}_j}, \qquad R := \sum_{j=1}^{n} z_j \frac{\partial}{\partial z_j} \tag{6}$$

the anti-holomorphic and the holomorphic radial derivative, respectively. Similarly, the product of Hankel operators $T = H^*_{g_1} H_{f_1} \dots H^*_{g_n} H_{f_n} = (T_{\overline{g}_1 f_1} - T_{\overline{g}_1} T_{f_1}) \dots (T_{\overline{g}_n f_n} - T_{\overline{g}_n} T_{f_n})$ belongs to $\mathcal{S}^{\mathrm{Dixm}}$ and

$$\mathrm{tr}_\omega T = \frac{1}{n!} \int_{\mathbf{S}^{2n-1}} \prod_{j=1}^{n} \{f_j, g_j\}_L \, d\sigma, \tag{7}$$

with the "half-bracket"

$$\{f, g\}_L = \sum_{j=1}^{n} \frac{\partial f}{\partial \overline{z}_j} \frac{\partial \overline{g}}{\partial z_j} - \overline{R} f R \overline{g}. \tag{8}$$

Note that for $n = 1$ the right-hand sides of (4) and (7) vanish, in accordance with the fact that the cutoff then occurs at $p = 1$ instead of $p = 2n = 2$; in fact, it was shown by Rochberg and one of the present authors that for $n = 1$ actually $|H_f| = (H^*_f H_f)^{1/2}$, rather than $H^*_f H_f$, is in the Dixmier class for any $f \in C^\infty(\overline{\mathbf{D}})$, and

$$\mathrm{tr}_\omega |H_f| = \int_{\partial \mathbf{D}} |\overline{\partial} f| \, d\sigma.$$

The result generalizes further from $\mathbf{B}^n$ to arbitrary bounded strictly pseudoconvex domains $\Omega \subset \mathbf{C}^n$ with smooth boundary, with the "half-bracket" $\{f, g\}_L$ and "boundary Poisson bracket" $\{f, g\}_b$ replaced by $\mathcal{L}(\overline{\partial}_b f, \overline{\partial}_b g)$ and $\mathcal{L}(\overline{\partial}_b f, \overline{\partial}_b \overline{g}) - \mathcal{L}(\overline{\partial}_b g, \overline{\partial}_b \overline{f})$, respectively, where $\mathcal{L}$ is the dual Levi form; see [12].

The aim of the present paper is to establish analogous results also in the context of the Fock space $\mathcal{F}_\gamma$. The role of functions smooth on the closure is now played by functions $f \in C^\infty(\mathbf{C}^n)$ admitting the asymptotics

$$f(z) \sim f_0\Big(\frac{z}{|z|}\Big) + |z|^{-1} f_1\Big(\frac{z}{|z|}\Big) + |z|^{-2} f_2\Big(\frac{z}{|z|}\Big) + \dots \tag{9}$$

as $|z| \to +\infty$, for some $f_0, f_1, \dots \in C^\infty(\mathbf{S}^{2n-1})$; see Section 3 below for the precise statement. Our main result on Hankel operators is then the following.

**Theorem 1.** *Let $f, g$ be as in* (9)*. Then $(H^*_f H_g)^n$ belongs to the Dixmier class, and*

$$\mathrm{tr}_\omega (H^*_f H_g)^n = \frac{1}{n!} \int_{\mathbf{S}^{2n-1}} \{\overline{f}_0, g_0\}^n_F \, d\sigma,$$

*where for $\phi, \psi \in C^\infty(\mathbf{S}^{2n-1})$,*

$$\{\phi, \psi\}_F = \sum_{j=1}^{n} \partial^b_j \phi \cdot \overline{\partial}^b_j \psi - E\phi \cdot E\psi.$$

*Here $\partial_j^b, \overline{\partial}_j^b$ and $E$ are vector fields on $\mathbf{S}^{2n-1}$ defined in Section 6.*

One also has an analogous result for $\operatorname{tr}_\omega H_f^* H_g$ if the first $n$ terms in the expansions (9) for $f, g$ identically vanish; or for $\operatorname{tr}_\omega T_f$ if $f_0 = f_1 = \dots = f_{2n-1} \equiv 0$ in (9); as well as for the Weyl operator $W_f$ with such $f$. See again Sections 5 and 6 below for the details.

The proof for the ball in [11] made use of the "pseudo-Toeplitz operators" of Howe [18], thus in effect — paradoxically — passing from the ball to the Fock space; thus, in principle, the same method could be used to handle also our case of $\mathcal{F}_\gamma$ here. However, we proceed by a simpler avenue: namely, using the standard correspondence between Toeplitz and Weyl calculi [13], we reduce the problem to one about the Weyl operators $W_f$, and then handle the latter by invoking the existing theory for the Weyl pseudodifferential operators available e.g. in [17], [6], [15], and especially [21] and [22].

The prerequisites on the Weyl calculus are reviewed in Section 2. Its relationship to Toeplitz operators is discussed in Section 3. Some Schatten-class properties of Weyl operators are handled in Section 4. The proof of Theorem 1 appears in Section 6, after some preparations in Section 5. Various concluding comments and open problems are collected in Section 7. Throughout the paper, we will be using the standard multi-index notations like $z^\alpha = z_1^{\alpha_1} \dots z_n^{\alpha_n}$, $\alpha! = \alpha_1! \dots \alpha_n!$, $\partial_z^\alpha = \partial^{|\alpha|}/\partial z_1^{\alpha_1} \dots \partial z_n^{\alpha_n}$, etc., for $\alpha = (\alpha_1, \dots, \alpha_n) \in \mathbf{N}^n$, $\mathbf{N} = \{0, 1, 2, \dots\}$. Furthermore, unless explicitly stated otherwise, the letters $z, w$ will always denote elements of $\mathbf{C}^n$, while $x, y, \xi$ are reserved for variables in $\mathbf{R}^n$; and, similarly, $\partial_x, \partial_y, \partial_\xi$ stand for ordinary partial derivatives, while the symbols $\partial_z$ and $\overline{\partial}_z$ denote the Wirtinger operators (holomorphic and anti-holomorphic derivatives).

## 2. The Weyl calculus

For a function $a(x, \xi)$ of $(x, \xi) \in \mathbf{R}^n \times \mathbf{R}^n$, the Weyl operator $W_a^{(\epsilon)} \equiv W_a$ on $L^2(\mathbf{R}^n)$ with symbol $a$ is defined by

$$W_a f(x) = \Big(\frac{|\epsilon|}{2\pi}\Big)^n \int_{\mathbf{R}^n} \int_{\mathbf{R}^n} a\Big(\frac{x+y}{2}, \xi\Big) \, e^{\epsilon i \langle x-y, \xi \rangle} \, f(y) \, dy \, d\xi. \tag{10}$$

Here $\epsilon \neq 0$ is a (real) parameter; the standard choices are $\epsilon = 2\pi$ [13], $\epsilon = -1$ [3] or $\epsilon = 1$ [22] [21]. From the Fourier inversion formula, one sees, in particular, that $W_a$ is just the operator of "multiplication by $a$" if $a(x, \xi) = a(x)$ depends only on $x$, while

$$W_{\xi^\alpha} = \frac{1}{(\epsilon i)^\alpha} \partial_x^\alpha.$$

The simple relation

$$W_a^{(\epsilon)} = W_{a_{\delta/\epsilon}}^{(\delta)}, \tag{11}$$

where $a_t(x, \xi) := a(x, t\xi)$, makes it possible to translate results for one value of $\epsilon$ into those for another. The symbol $a$ in (10) has to be a sufficiently nice function (e.g. Schwartz) in order for the integral to make sense; however, standard arguments show that the definition can be extended to all tempered distributions $a$ on $\mathbf{R}^{2n}$, and any continuous operator from $\mathcal{S}(\mathbf{R}^n)$ into $\mathcal{S}'(\mathbf{R}^n)$ is then of the form $W_a$ for a unique $a \in \mathcal{S}'(\mathbf{R}^{2n})$ [13]. To get nice properties, one restricts $a$ to lie in various symbol classes, involving decay conditions on derivatives of $a$ as either $|\xi|$, or both $|x|$ and $|\xi|$, tend to infinity [22] [15] [6] [17] [21]; for our purpose, the convenient symbol classes are those of Shubin [21] and Voros [22]. Namely, for $m \in \mathbf{R}$, define

$$GLS^m := \{a \in C^\infty(\mathbf{R}^n \times \mathbf{R}^n) : \|a\|_{\alpha\beta m} < \infty \,\forall \alpha, \beta\},$$

where, for multiindices $\alpha, \beta \in \mathbf{N}^n$,

$$\|a\|_{\alpha\beta m} = \sup_{x,\xi} |\partial_x^\alpha \partial_\xi^\beta a(x,\xi)|(1 + |x|^2 + |\xi|^2)^{(|\alpha|+|\beta|-m)/2}.$$

The spaces $GLS^m$ clearly increase with $m$, and we denote $GLS^\infty := \bigcup_{m\in\mathbf{R}} GLS^m$ and $GLS^{-\infty} := \bigcap_{m\in\mathbf{R}} GLS^m = \mathcal{S}(\mathbf{R}^{2n})$, the Schwartz space on $\mathbf{R}^{2n}$. As usual, one says that $a \in GLS^m$ has an asymptotic expansion $a \sim \sum_j a_j$ if $a_j \in GLS^{m_j}$, $m_j \searrow -\infty$, $m_0 = m$ and

$$a - \sum_{j=0}^{N-1} a_j \in GLS^{m_N}, \qquad \forall N = 0, 1, 2, \dots, \tag{12}$$

and shows that for any sequence $m_j \searrow -\infty$ and $a_j \in GLS^{m_j}$, $j \in \mathbf{N}$, there exists $a \in GLS^{m_0}$, unique modulo $GLS^{-\infty}$, for which (12) holds.

For $a, b \in GLS^\infty$, the operators $W_a, W_b$ can be shown to map $\mathcal{S}$ continuously into itself, hence so does $W_a W_b$ and thus $W_a W_b = W_c$ for a unique $c := a\#b \in \mathcal{S}'(\mathbf{R}^{2n})$; and one has $a\#b \in GLS^{m+k}$ if $a \in GLS^m$, $b \in GLS^k$, and

$$a\#b \sim \sum_{j=0}^\infty \frac{(i/2)^j}{j!\epsilon^j} (\partial_{x,a}\partial_{\xi,b} - \partial_{\xi,a}\partial_{x,b})^j ab, \tag{13}$$

where the subscript $a$ in $\partial_{\xi,a}$ means that $\partial_\xi$ is applied only to $a$, and analogously for $\partial_{x,b}, \partial_{x,a}$ and $\partial_{\xi,b}$. (For $\epsilon = 1$, (13) is Theorem 2.4.1 in [22], and for $\epsilon = 2\pi$, it is the formula (2.48) in [13]; the result for arbitrary $\epsilon$ then follows by (11).) Note that in this case (12) holds with $m_j = m + k - 2j$.

One can also define Weyl operators on the Fock spaces $\mathcal{F}_\gamma$. Namely, for $\alpha \in \mathbf{R} \setminus \{0\}$, the *Bargmann transform*

$$B_{\alpha\gamma} f(z) \equiv Bf(z) := c_{\alpha\gamma} \int_{\mathbf{R}^n} f(x) e^{\alpha x z - \frac{\alpha^2}{4\gamma} x^2 - \frac{\gamma}{2} z^2} \, dx,$$

where

$$c_{\alpha\gamma} = \Big(\frac{\alpha^2\gamma}{2\pi^3}\Big)^{n/4},$$

defines an isometry of $L^2(\mathbf{R}^n)$ onto $\mathcal{F}_\gamma$, with inverse

$$B^{-1}F(x) = c_{\alpha\gamma} \int_{\mathbf{C}^n} F(z) e^{\alpha x \bar z - \frac{\alpha^2}{4\gamma} x^2 - \frac{\gamma}{2} \bar z^2} \, d\mu_\gamma(z).$$

(Here $z^2 := \sum_{j=1}^n z_j^2 = \langle z, \bar z\rangle$ and similarly $x^2 = \langle x, x\rangle$, $x\bar z = \langle x, z\rangle$ and $xz = \langle x, \bar z\rangle$.) For $\alpha = 2\pi$, $\gamma = \pi$ this is established e.g. in Chapter 1, §6 of [13]; for arbitrary $\gamma > 0$, $\alpha \in \mathbf{R} \setminus \{0\}$, this again follows from the simple relations

$$\delta_t B_{\alpha\gamma} = s^n B_{\alpha s t, \gamma t^2} \delta_s, \qquad s, t \in \mathbf{R} \setminus \{0\}, \quad \delta_t f(x) := f(tx),$$

and a simple change-of-variable argument.

*Remark.* In the literature, one often encounters Fock spaces $\widetilde{\mathcal{F}_\gamma}$ with respect to the normalized Gaussians (of total mass one) $d\tilde\mu_\gamma(z) = (\frac{\gamma}{\pi})^n d\mu_\gamma(z)$. The constant $c_{\alpha\gamma}$ then must be replaced by $\tilde c_{\alpha\gamma} = \Big(\dfrac{\alpha^2}{2\pi\gamma}\Big)^{n/4}$. □

A simple calculation shows that

$$(B^{-1}z_jBf)(x) = \Big(\frac{\alpha}{2\gamma}x_j - \frac{1}{\alpha}\frac{\partial}{\partial x_j}\Big)f(x),$$

or, since the Toeplitz operator $T_{z_j}$ is just the multiplication by $z_j$ on $\mathcal{F}_\gamma$, $B^{-1}T_{z_j}B = \frac{\alpha}{2\gamma}x_j - \frac{1}{\alpha}\frac{\partial}{\partial x_j}$. Similarly $T_{\overline{z}_j} = \frac{1}{\gamma}\frac{\partial}{\partial z_j}$ on $\mathcal{F}_\gamma$ and

$$B^{-1}T_{\overline{z}_j}B = \frac{\alpha}{2\gamma}x_j + \frac{1}{\alpha}\frac{\partial}{\partial x_j}.$$

Thus if $\epsilon < 0$ and we set

$$\gamma = -\frac{\epsilon}{2}, \qquad \alpha = -\epsilon,$$

then $B = B_{-\epsilon,-\epsilon/2}$ will map $L^2(\mathbf{R}^n)$ unitarily onto $\mathcal{F}_\gamma(\mathbf{C}^n)$ and satisfy

$$\begin{aligned} B^{-1}T_{z_j}B &= x_j + \frac{1}{\epsilon}\frac{\partial}{\partial x_j} = W_{x_j+i\xi_j}, \\ B^{-1}T_{\overline{z}_j}B &= x_j - \frac{1}{\epsilon}\frac{\partial}{\partial x_j} = W_{x_j-i\xi_j}. \end{aligned} \tag{14}$$

We are thus led to define, with $B = B_{2\gamma,\gamma}$, the Weyl operators $W_a^{(\gamma)} \equiv W_a$ on $\mathcal{F}_\gamma$ by

$$W_a^{(\gamma)} := BW^{(\epsilon)}B^{-1}, \qquad \epsilon = -2\gamma,$$

where the $W_a^{(\epsilon)}$ on the right-hand side is the Weyl operator (10) on $L^2(\mathbf{R}^n)$. From now on, we will also view $a(x,\xi)$ rather as the function $a(z)$ of $z = x + i\xi \in \mathbf{C}^n$. The relations (14) then become simply

$$W_{z_j}^{(\gamma)} = T_{z_j}, \quad W_{\overline{z}_j}^{(\gamma)} = T_{\overline{z}_j}.$$

Note that the symbol classes $GLS^m$ can now be written as

$$GLS^m = \{a \in C^\infty(\mathbf{C}^n) : \sup_z |\partial^\alpha\overline{\partial}^\beta a(z)|(1+|z|^2)^{(|\alpha|+|\beta|-m)/2} < \infty,\ \forall \alpha,\beta \in \mathbf{N}^n\}.$$

Similarly, in view of the relations $\partial_x = \partial + \overline{\partial}$, $\partial_\xi = i(\partial - \overline{\partial})$, (13) becomes

$$\begin{aligned} a\#b &\sim \sum_{j=0}^\infty \frac{1}{j!\epsilon^j}(\partial_a\overline{\partial}_b - \partial_b\overline{\partial}_a)^j ab \\ &= \sum_{\alpha,\beta} \frac{(-1)^{|\beta|}}{\alpha!\beta!(-2\gamma)^{|\alpha|+|\beta|}}\partial^\alpha\overline{\partial}^\beta a \cdot \partial^\beta\overline{\partial}^\alpha b, \end{aligned} \tag{15}$$

where the subscripts at $\partial, \overline{\partial}$ on the first line indicate to which function it only applies. Note that for $a \in GLS^m$, $b \in GLS^k$, the general term in the last sum belongs to $GLS^{m+k-2|\alpha|-2|\beta|}$.

*Remark.* For $\epsilon < 0$, one instead of (14) obtains that $B = B_{\epsilon,\epsilon/2}$ satisfies $B^{-1}T_{z_j}B = W_{\overline{z}_j}$ and $B^{-1}T_{\overline{z}_j}B = W_{z_j}$; this is the situation for $\epsilon = 2\pi$ in [13]. □

## 3. Transition between Toeplitz and Weyl operators

Let $\mathcal{E} = e^{\Delta/8\gamma}$ denote the heat solution operator on $\mathbf{C}^n$ at time $t = \frac{1}{8\gamma}$:

$$\mathcal{E}f(z) = \Big(\frac{2\gamma}{\pi}\Big)^n \int_{\mathbf{C}^n} f(w) e^{-2\gamma|z-w|^2} \, dw.$$

Furthermore, denoting by $k_w(z) = (\frac{\gamma}{\pi})^{n/2} e^{\gamma\langle z,w\rangle - \gamma|w|^2/2}$ the normalized reproducing kernel of $\mathcal{F}_\gamma$ at $w \in \mathbf{C}^n$, it is known that any bounded linear operator $T$ on $\mathcal{F}_\gamma$ is uniquely determined by its *Berezin symbol*

$$\widetilde{T}(w) = \langle Tk_w, k_w\rangle_{\mathcal{F}_\gamma} = e^{-\gamma|w|^2} (Te^{\gamma\langle\,\cdot\,,w\rangle})(w).$$

For Toeplitz operators, a simple calculation gives

$$\widetilde{T_f}(z) = \Big(\frac{\gamma}{\pi}\Big)^n \int_{\mathbf{C}^n} f(w) e^{-\gamma|z-w|^2} \, dw = \mathcal{E}^2 f(z),$$

while a slightly more difficult argument gives

$$\widetilde{W_f} = \mathcal{E}f.$$

Comparing the last two formulas and appealing to the uniqueness, it follows that

$$T_f = W_{\mathcal{E}f}. \tag{16}$$

(For $\gamma = \pi$, the proof of all this can be found in Chapter 2, §9 in [13]; the case of general $\gamma > 0$ follows again by the change of variable $f(z) \mapsto f(tz)$, $t \neq 0$, which maps $\mathcal{F}_\gamma$ unitarily — up to a constant factor — onto $\mathcal{F}_{\gamma/t^2}$.)

The Taylor series $e^{\Delta/8\gamma} = I + \frac{\Delta}{8\gamma} + \frac{\Delta^2}{2!8^2\gamma^2} + \dots$ suggests that for $a \in GLS^m$,

$$\mathcal{E}a \sim \sum_{j=0}^{\infty} \frac{\Delta^j a}{j!(8\gamma)^j},$$

with the $j$-th term on the right belonging to $GLS^{m-2j}$. While this may be true in full generality, we will need this only for a special class of functions. Let $\mathbf{S}^{2n-1}$ denote the unit sphere in $\mathbf{C}^n \cong \mathbf{R}^{2n}$.

**Definition.** For $m \in \mathbf{R}$, let $\mathcal{A}^m$ be the space of all functions in $C^\infty(\mathbf{C}^n)$ for which

$$f(z) \sim \sum_{j=0}^{\infty} |z|^{m-j} f_j\Big(\frac{z}{|z|}\Big) \qquad \text{as } |z| \to +\infty, \tag{17}$$

for some $f_j \in C^\infty(\mathbf{S}^{2n-1})$, $j \in \mathbf{N}$, in the sense that

$$\sup_{|z|\geq 1} |z|^{N+|\alpha|+|\beta|-m} \partial^\alpha \overline{\partial}^\beta \Big[ f(z) - \sum_{j=0}^{N-1} |z|^{m-j} f_j\Big(\frac{z}{|z|}\Big)\Big] < \infty \tag{18}$$

for any multiindices $\alpha, \beta$.

In other words, $f \in \mathcal{A}^m$ if $f \in GLS^m$ and has an asymptotic expansion

$$f \sim \sum_{j=0}^{\infty} a_j \tag{19}$$

where $a_j \in GLS^{m-j}$ is homogeneous of degree $m - j$ for $|z| \geq 1$.

An example of a function in $\mathcal{A}^m$ is, of course, any smooth function on $\mathbf{C}^n$ equal to $f_m(\frac{z}{|z|})|z|^m$ for $|z| \geq 1$, with some $f_m \in C^\infty(\mathbf{S}^{2n-1})$.

**Theorem 2.** *For $f \in \mathcal{A}^m$, $m \le 0$,*

$$\mathcal{E}f \sim \sum_{k=0}^{\infty} \sum_{\substack{j,l\ge 0\\ j+2l=k}} \frac{\Delta^l}{l!(8\gamma)^l} |z|^{m-j} f_j\Big(\frac{z}{|z|}\Big), \tag{20}$$

*where $f_j$ are the functions from* (17).

Note that if a function is homogeneous of degree $d$, then its derivatives of order $l$ are homogeneous of degree $d-l$. Thus it is immediate that the $k$-th term in the sum on the right is homogeneous of degree $m-k$ for $|z|\ge 1$, and, hence, indeed, belongs to $GLS^{m-k}$.

The proof relies on the following two lemmas, which in fact hold for any $m \in \mathbf{R}$, but the simpler versions below are sufficient for our purposes.

**Lemma 3.** *Let $f \in L^\infty(\mathbf{C}^n)$ satisfy $|f(z)| \le C_1 \cdot (1+|z|^2)^{m/2}$, with some $m \le 0$, and let $g$ be a Schwartz function. Then there is $C_2 > 0$ such that the convolution $f * g$ satisfies $|(f*g)(z)| \le C_2 \cdot (1+|z|^2)^{m/2}$.*

*Proof.* For brevity, denote $b = -\frac{m}{2} \ge 0$, and choose some $a > b+n$. Since $g \in \mathcal{S}$, there is $C_3 > 0$ such that $|g(z)| \le C_3 \cdot (1+|z|^2)^{-a}$. From the elementary inequality

$$1+|x+y|^2 \le 2(1+|x|^2)(1+|y|^2), \qquad x,y \in \mathbf{C}^n,$$

we have

$$1+|w|^2 \le 2(1+|z|^2)(1+|z-w|^2),$$

or

$$\frac{1}{1+|z-w|^2} \le \frac{2(1+|z|^2)}{1+|w|^2}, \qquad z,w \in \mathbf{C}^n.$$

Consequently,

$$\begin{aligned} |(f*g)(w)| &\le \int_{\mathbf{C}^n} |f(w-z)g(z)|\,dz \\ &\le C_1C_3 \int_{\mathbf{C}^n} (1+|z|^2)^{-a}(1+|z-w|^2)^{-b}\,dz \\ &\le \frac{2^bC_1C_3}{(1+|w|^2)^b} \int_{\mathbf{C}^n} (1+|z|^2)^{b-a}\,dz. \end{aligned}$$

Since $b-a < -n$, the last integral is finite, and the assertion follows. $\square$

**Lemma 4.** *Let $f \in L^\infty(\mathbf{C}^n)$ satisfy $|f(z)| \le C_1 \cdot (1+|z|^2)^{m/2}$, with some $m \le 0$, and let $g$ be supported on $|z| \le 1$ and integrable. Then there is $C_2 > 0$ such that the convolution $f*g$ satisfies $|(f*g)(z)| \le C_2 \cdot (1+|z|^2)^{m/2}$.*

*Proof.* Immediate from

$$\begin{aligned} |(f*g)(w)| &\le \int_{\mathbf{C}^n} |f(w-z)g(z)|\,dz \\ &= C_1 \int_{|z|<1} (1+|z-w|^2)^{m/2}|g(z)|\,dz \\ &\le C_1\|g\|_1 \sup_{|z|<1} (1+|z-w|^2)^{m/2} \\ &\le C_2(1+|w|^2)^{m/2}. \quad \square \end{aligned}$$

*Proof of Theorem 2.* Note that $\mathcal{E}$ is the operator of convolution with the Schwartzian function $(\frac{2\gamma}{\pi})^n e^{-2\gamma|z|^2}$; by standard properties of convolutions [20, Theorem 7.19], we thus have for any multiindices $\alpha, \beta$

$$\mathcal{E}(\partial^\alpha \overline{\partial}^\beta f) = \partial^\alpha \overline{\partial}^\beta (\mathcal{E} f).$$

Furthermore, $\mathcal{E}$ is an integral operator with positive kernel, so $|f| \leq g$ implies $|\mathcal{E} f| \leq \mathcal{E} g$.

Choose $a_j \in C^\infty(\mathbf{C}^n)$ such that $a_j(z) = |z|^{m-j} f_j(\frac{z}{|z|})$ for $|z| \geq 1$ (so that (19) holds). Then for any $N = 0, 1, 2, \dots,$

$$\begin{aligned} \partial^\alpha \overline{\partial}^\beta \Big[\mathcal{E} f - \sum_{j=0}^{N-1} \mathcal{E} a_j\Big] &= \mathcal{E}\Big(\partial^\alpha \overline{\partial}^\beta \Big[f - \sum_{j=0}^{N-1} a_j\Big]\Big) \\ &= \mathcal{E}(O((1+|z|^2)^{(m-|\alpha|-|\beta|-N)/2}) \\ &= O((1+|z|^2)^{(m-|\alpha|-|\beta|-N)/2}) \end{aligned}$$

by Lemma 3. Thus $\mathcal{E} f - \sum_{j=0}^{N-1} \mathcal{E} a_j \in GLS^{m-N}$ and

$$\mathcal{E} f \sim \sum_{j=0}^{\infty} \mathcal{E} a_j.$$

It is therefore enough to prove the theorem for each $a_j$ separately; that is, we may from now on assume that $f \in C^\infty(\mathbf{C}^n)$ satisfies $f(z) = |z|^m f_0(\frac{z}{|z|})$ for $|z| \geq 1$, with some $m \leq 0$ and $f_0 \in C^\infty(\mathbf{S}^{2n-1})$. Using the usual notations $^\wedge$ and $^\vee$ for the Fourier transform and its inverse, respectively, we then have

$$\begin{aligned} \mathcal{E} f - \sum_{j=0}^{N-1} \frac{\Delta^j f}{j!(8\gamma)^j} &= \Big[\widehat{f} \cdot \Big(e^{-|\xi|^2/8\gamma} - \sum_{j=0}^{N-1} \frac{(-|\xi|^2)^j}{j!(8\gamma)^j}\Big)\Big]^\vee \\ &= [(-|\xi|^2)^N \widehat{f} \cdot \widehat{s}]^\vee \\ &= \Delta^N f * s, \end{aligned}$$

where

$$s = \left[\frac{e^{-|\xi|^2/8\gamma} - \sum_{j=0}^{N-1} \frac{(-|\xi|^2)^j}{j!(8\gamma)^j}}{(-|\xi|^2)^N}\right]^\vee.$$

Clearly, $\widehat{s} = \psi + \sum_{j=1}^{N} \psi_j$, where $\psi \in \mathcal{S}$ while $\psi_j(\xi) = 0$ for $|\xi| \leq \frac{1}{2}$ and

$$\psi_j(\xi) = \frac{(-1)^j}{(N-j)!(8\gamma)^{(N-j)}} |\xi|^{-2j} \qquad \text{for } |\xi| \geq 1.$$

It is known (see e.g. [19], Proposition 2.1) that the inverse Fourier transform of $\psi_j$ has the form

$$\psi_j^\vee(z) = a_j(z) + \chi(|z|) \cdot \begin{cases} c_j |z|^{2j-2n} & \text{if } 2j-2n < 0 \\ c_j |z|^{2j-2n} \log|z| & \text{if } 2j-2n \geq 0 \end{cases},$$

where $\chi(|z|)$ is some fixed smooth cutoff function vanishing on $|z| \geq 1$ and equal to 1 on $|z| \leq \frac{1}{2}$, and $a_j \in \mathcal{S}$. Altogether, we see that

$$s = s_1 + s_2,$$

where $s_1$ is a Schwartz function, while $s_2$ is supported on $|z| < 1$ and integrable. Replacing $f$ by $\partial^\alpha \overline{\partial}^\beta f$ and using again Lemma 3 (for convolution with $s_1$) as well as Lemma 4 (for convolution with $s_2$), we get

$$\begin{aligned}
\partial^\alpha \overline{\partial}^\beta \Big[\mathcal{E} f - \sum_{j=0}^{N-1} \frac{\Delta^j f}{j!(8\gamma)^j}\Big] &= \Big[\mathcal{E} - \sum_{j=0}^{N-1} \frac{\Delta^j}{j!(8\gamma)^j}\Big](\partial^\alpha \overline{\partial}^\beta f) \\
&= (\Delta^N \partial^\alpha \overline{\partial}^\beta f) * s \\
&= O((1+|z|^2)^{(m-|\alpha|-|\beta|-2N)/2}) * s \\
&= O((1+|z|^2)^{(m-|\alpha|-|\beta|-2N)/2}).
\end{aligned}$$

Thus $\mathcal{E} f - \sum_{j=0}^{N-1} \frac{\Delta^j f}{j!(8\gamma)^j} \in GLS^{m-2N}$, and

$$\mathcal{E} f \sim \sum_{j=0}^{\infty} \frac{\Delta^j f}{j!(8\gamma)^j},$$

completing the proof of the theorem. □

## 4. Schatten and Dixmier classes

Recall that for a compact (linear) operator $A$ from one Hilbert space into another, the $s$-numbers $s_0(A) \geq s_1(A) \geq s_2(A) \geq \dots$ of $A$ are defined to be the eigenvalues of $|A| = (A^*A)^{1/2}$ arranged in non-increasing order (counting multiplicities). One has $s_0(A) = \|A\|$, the operator norm of $A$, and $s_j(A) \searrow 0$ as $j \to \infty$; also

$$s_j(A^*) = s_j(A), \quad s_j(A+B) \leq s_j(A) + s_j(B),$$

and

$$s_{j+k}(AB) \leq s_j(A) s_k(B); \tag{21}$$

in particular,

$$s_j(AB) \leq s_j(A)\|B\|, \qquad s_j(AB) \leq \|A\| s_j(B). \tag{22}$$

The Schatten classes $\mathcal{S}^p$, $0 < p < \infty$, are defined to consist of all $A$ with

$$\sum_{j=0}^{\infty} s_j(A)^p < \infty.$$

For $p \geq 1$, the $p$-th root of the left-hand side defines a norm on $\mathcal{S}^p$ which makes it into a Banach space (for $0 < p < 1$, one gets quasi-Banach spaces). For $p > 1$, one also has the Schatten classes $\mathcal{S}^{p,\infty}$ consisting of all $A$ with

$$\sup_j (j+1)^{1/p} s_j(A) < \infty,$$

with the left-hand side again making $\mathcal{S}^{p,\infty}$ into a Banach space. For $p = 1$, the left-hand side defines a norm which however is not complete; to get a Banach space, one needs to define $\mathcal{S}^{1,\infty} \equiv \mathcal{S}^{\mathrm{Dixm}}$ — the *Dixmier class* — as the space of all $A$ for which

$$\sup_{k\geq 2} \frac{\sum_{j=0}^{k} s_j(A)}{\log k} < \infty$$

with the corresponding norm. It follows from (22) that all the spaces $\mathcal{S}^p$, $0 < p < \infty$, and $\mathcal{S}^{p,\infty}$, $1 \leq p < \infty$, are ideals, i.e. if $B$ is bounded and $A \in \mathcal{S}^p$ then $AB, BA \in \mathcal{S}^p$ with $\|AB\|_p \leq \|A\|_p\|B\|$ and $\|BA\|_p \leq \|A\|_p\|B\|$, and similarly for $\mathcal{S}^{p,\infty}$. Likewise, from (21) one sees that

$$A, B \in \mathcal{S}^{2,\infty} \implies AB \in \mathcal{S}^{1,\infty}. \tag{23}$$

A good reference for all the material above is [14].

For operators in $\mathcal{S}^{\mathrm{Dixm}}$, one can define the so-called *Dixmier trace* $\operatorname{tr}_\omega$. Its construction involves a choice of a Banach limit $\omega$ on the space $\ell^\infty$ of all bounded sequences $\{a_j\}_{j=0}^\infty$, which is "scaling-invariant" in the sense that $\omega(a_1, a_2, \dots) = \omega(a_1, a_1, a_2, a_2, \dots)$. One then defines $\operatorname{tr}_\omega$ for positive operators in $\mathcal{S}^{\mathrm{Dixm}}$ by

$$\operatorname{tr}_\omega A = \omega\Big(\frac{1}{\log(k+2)} \sum_{j=0}^{k} s_j(A)\Big).$$

The scaling-invariance of $\omega$ turns out to guarantee that $\operatorname{tr}_\omega(A+B) = \operatorname{tr}_\omega A + \operatorname{tr}_\omega B$ for $A, B \geq 0$, and thus one can extend $\operatorname{tr}_\omega$ unambiguously to all $A \in \mathcal{S}^{\mathrm{Dixm}}$ by linearity. The resulting functional satisfies

$$\begin{aligned} \operatorname{tr}_\omega(A+B) &= \operatorname{tr}_\omega A + \operatorname{tr}_\omega B, \qquad \forall A, B \in \mathcal{S}^{\mathrm{Dixm}}, \\ \operatorname{tr}_\omega(AB) &= \operatorname{tr}_\omega(BA), \qquad \forall A \in \mathcal{S}^{\mathrm{Dixm}}, B \text{ bounded}, \end{aligned} \tag{24}$$

and

$$|\operatorname{tr}_\omega(AB)| \leq \|A\|_{1,\infty}\|B\| \qquad \forall A \in \mathcal{S}^{\mathrm{Dixm}}, B \text{ bounded}. \tag{25}$$

The value of $\operatorname{tr}_\omega A$ will in general depend on the choice of $\omega$; operators $A$ for which all choices of $\omega$ yield the same value are called *measurable*. Clearly, this includes, in particular, all operators $A \geq 0$ for which the limit

$$\lim_{k\to\infty} \frac{1}{\log k} \sum_{j=0}^{k} s_j(A) \tag{26}$$

exists, and $\operatorname{tr}_\omega A$ is then equal to this limit. Note that it follows from the construction that, in particular, $\operatorname{tr}_\omega A = 0$ whenever $A$ is trace class.

For more details and information about the Dixmier trace, see the book [8].

It has been shown by Voros [22, Theorem 2.7.1] (cf. [21], §27) that

$$W_a \in \mathcal{S}^2 \quad \text{if } a \in GLS^m, m < -n, \tag{27}$$

and

$$W_a \in \mathcal{S}^1 \quad \text{if } a \in GLS^m, m < -2n. \tag{28}$$

Using Theorem 2, this easily gives also a criterion for membership in the Dixmier class.

**Theorem 5.** *Let $a \in \mathcal{A}^{-2n}$. Then $W_a, T_a \in \mathcal{S}^{\rm Dixm}$. Also, $\operatorname{tr}_\omega T_a = \operatorname{tr}_\omega W_a$ for any choice of $\omega$.*

*Proof.* By (16) and Theorem 2,

$$T_a = W_{\mathcal{E}a} = W_a + W_b,$$

where $b \in GLS^{-2n-2}$. By (28), $W_b$ is trace-class, so $W_b \in \mathcal{S}^{\rm Dixm}$ and $\operatorname{tr}_\omega W_b = 0$. It thus suffices to show that $T_a \in \mathcal{S}^{\rm Dixm}$. We will even show that, in fact, $T_a \in \mathcal{S}^{\rm Dixm}$ whenever $a(z) = O(|z|^{-2n})$ as $|z| \to +\infty$ (i.e. one needs the estimate (18) from the definition of $\mathcal{A}^{-2n}$ only for $\alpha = \beta = N = 0$).

Indeed, one then has

$$a(z) = g(z)(1+|z|^2)^{-n}$$

with some $g \in L^\infty(\mathbf{C}^n)$. Now the operator $T_a$ can be identified, upon extension by zero (which has no effect on the $s$-numbers), with the operator $P_\gamma M_a P_\gamma$ on all of $L^2(\mathbf{C}^n, d\mu_\gamma)$, where $M_a : f \mapsto af$ stands for the operator of "multiplication by $a$". Now

$$P_\gamma M_a P_\gamma = P_\gamma M_{(1+|z|^2)^{-n/2}} M_g M_{(1+|z|^2)^{-n/2}} P_\gamma.$$

Since $M_g$ is bounded, the fact that $\mathcal{S}^{2,\infty}$ is an ideal and (23) imply that $T_a$ belongs to $\mathcal{S}^{1,\infty}$ if $M_{(1+|z|^2)^{-n/2}} P_\gamma$ belongs to $\mathcal{S}^{2,\infty}$. Since $(M_{(1+|z|^2)^{-n/2}} P_\gamma)^* M_{(1+|z|^2)^{-n/2}} P_\gamma = T_{(1+|z|^2)^{-n}}$, we need to check that the eigenvalues $\{c_j\}_{j=0}^\infty$ of $T_{(1+|z|^2)^{-n}}$, arranged in non-increasing order (counting multiplicities), satisfy $\sup_j (j+1)^{1/2}\sqrt{c_j} < \infty$, or

$$\sup_j (j+1) c_j < \infty.$$

However, it is well known that the monomials $\{z^\alpha\}_{\alpha \text{ a multiindex}}$ form an orthogonal basis of $\mathcal{F}_\gamma$, with norms

$$\|z^\alpha\|^2_{\mathcal{F}_\gamma} = \Big(\frac{\pi}{\gamma}\Big)^n \frac{\alpha!}{\gamma^{|\alpha|}}.$$

Furthermore, denoting by $d\zeta$ the surface measure on $\mathbf{S}^{2n-1}$,

$$\begin{aligned}
\langle T_{(1+|z|^2)^{-n}} z^\alpha, z^\beta \rangle &= \int_{\mathbf{C}^n} z^\alpha \overline{z}^\beta (1+|z|^2)^{-n} e^{-\gamma|z|^2}\, dz \\
&= \int_0^\infty r^{2|\alpha|}(1+r^2)^{-n} e^{-\gamma r^2} r^{2n-1}\, dr \int_{\mathbf{S}^{2n-1}} \zeta^\alpha \overline{\zeta}^\beta \, d\zeta \\
&= \frac{\int_0^\infty r^{2|\alpha|}(1+r^2)^{-n} e^{-\gamma r^2} r^{2n-1}\, dr}{\int_0^\infty r^{2|\alpha|} e^{-\gamma r^2} r^{2n-1}\, dr} \int_{\mathbf{C}^n} z^\alpha \overline{z}^\beta e^{-\gamma|z|^2}\, dz \\
&= \frac{\int_0^\infty t^{|\alpha|}(1+t)^{-n} e^{-\gamma t} t^{n-1}\, dt}{\int_0^\infty t^{|\alpha|} e^{-\gamma t} t^{n-1}\, dt} \langle z^\alpha, z^\beta \rangle_{\mathcal{F}_\gamma} \\
&= \int_0^\infty t^{|\alpha|+n-1}(1+t)^{-n} e^{-\gamma t}\, dt \frac{\gamma^{|\alpha|+n}}{\Gamma(|\alpha|+n)} \delta_{\alpha\beta} \|z^\alpha\|^2_{\mathcal{F}_\gamma}.
\end{aligned}$$

It follows that $T_{(1+|z|^2)^{-n}}$ is diagonalized by the monomial basis: $T_{(1+|z|^2)^{-n}} z^\alpha =$

$c_\alpha z^\alpha$, with

$$
\begin{aligned}
c_\alpha &= \frac{\gamma^{|\alpha|+n}}{\Gamma(|\alpha|+n)} \int_0^\infty \frac{t^{|\alpha|+n-1}}{(1+t)^n}\, e^{-\gamma t}\, dt \\
&\le \frac{\gamma^{|\alpha|+n}}{\Gamma(|\alpha|+n)} \int_0^\infty t^{|\alpha|-1}\, e^{-\gamma t}\, dt \\
&= \frac{\gamma^{|\alpha|+n}}{\Gamma(|\alpha|+n)} \frac{\Gamma(|\alpha|)}{\gamma^{|\alpha|}} = \frac{\gamma^n}{|\alpha|(|\alpha|+1)\dots(|\alpha|+n-1)} \\
&\le \Big(\frac{\gamma}{|\alpha|}\Big)^n
\end{aligned}
\tag{29}
$$

for $\alpha \neq 0$. Now the number of multiindices $\alpha$ with $|\alpha| = k$ is $\binom{k+n-1}{k}$; thus upon rearranging $c_\alpha$ in non-increasing order to get $c_j$, we get for $|\alpha| = k$

$$
\begin{aligned}
(j+1)c_j &\le \Big[1 + \sum_{l=0}^{k} \binom{l+n-1}{l}\Big] \Big(\frac{\gamma}{k}\Big)^n \\
&= \Big[1 + \binom{k+n}{n}\Big] \Big(\frac{\gamma}{k}\Big)^n \\
&\approx \frac{k^n}{n!}\Big(\frac{\gamma}{k}\Big)^n = \frac{\gamma^n}{n!} < \infty,
\end{aligned}
\tag{30}
$$

proving the result. □

For later use, we remark that the computations in (29) and (30) are, in some sense, sharp. Namely, for $t \ge \delta > 0$, we have $\frac{1}{t+1} \ge \frac{\delta/t}{\delta+1}$, so (writing $|\alpha| = k$ for brevity, $k > 0$)

$$
\begin{aligned}
\int_\delta^\infty \frac{t^{k+n-1}}{(1+t)^n}\, e^{-\gamma t}\, dt &\ge \frac{\delta^n}{(\delta+1)^n} \int_\delta^\infty t^{k-1} e^{-\gamma t}\, dt \\
&= \frac{\delta^n}{(\delta+1)^n} \Big(\frac{\Gamma(k)}{\gamma^k} - \int_0^\delta t^{k-1} e^{-\gamma t}\, dt\Big) \\
&\ge \frac{\delta^n}{(\delta+1)^n} \Big(\frac{\Gamma(k)}{\gamma^k} - \frac{\delta^k}{e^{\gamma\delta}}\Big)
\end{aligned}
$$

for $k > \gamma\delta + 1$, since the function $t^{k-1}e^{-\gamma t}$ is increasing on the interval $(0, \frac{k-1}{\gamma})$. Thus

$$
\begin{aligned}
c_\alpha &= \frac{\gamma^{k+n}}{\Gamma(k+n)} \int_0^\infty \frac{t^{k+n-1}}{(1+t)^n}\, e^{-\gamma t}\, dt \\
&\ge \frac{\gamma^{k+n}}{\Gamma(k+n)} \frac{\delta^n}{(\delta+1)^n} \Big(\frac{\Gamma(k)}{\gamma^k} - \frac{\delta^k}{e^{\gamma\delta}}\Big),
\end{aligned}
$$

for $k > \gamma\delta + 1$, and

$$
\liminf_{k\to\infty} \Big(\frac{k}{\gamma}\Big)^n c_\alpha \ge \liminf_{k\to\infty} \frac{\delta^n}{(\delta+1)^n} \Big(1 - \frac{\gamma^k \delta^k}{\Gamma(k) e^{\gamma\delta}}\Big) = \frac{\delta^n}{(\delta+1)^n}.
$$

On the other hand, by (29), $(\frac{k}{\gamma})^n c_\alpha \le 1$; letting $\delta \to +\infty$, we thus get

$$\lim_{k\to\infty} \Big(\frac{k}{\gamma}\Big)^n c_\alpha = 1.$$

Thus in (30), for $j$ large enough,

$$\begin{aligned}(j+1)c_j &\gtrsim \Big[1+\sum_{l=0}^{k}\binom{l+n-1}{l}\Big]\Big(\frac{\gamma}{k}\Big)^n \\ &= \Big[1+\binom{k+n}{n}\Big]\Big(\frac{\gamma}{k}\Big)^n \\ &\approx \frac{k^n}{n!}\Big(\frac{\gamma}{k}\Big)^n = \frac{\gamma^n}{n!},\end{aligned}$$

implying that

$$\lim_{j\to\infty}(j+1)c_j = \frac{\gamma^n}{n!}. \tag{31}$$

In particular,

$$\lim_{k\to\infty}\frac{1}{\log k}\sum_{j=0}^{k} c_j = \lim_{k\to\infty}\frac{1}{\log k}\frac{\gamma^n}{n!}\sum_{j=0}^{k}\frac{1}{j+1} = \frac{\gamma^n}{n!},$$

so by (26) the operator $T_{(1+|z|^2)^{-n}}$ is measurable and

$$\operatorname{tr}_\omega T_{(1+|z|^2)^{-n}} = \frac{\gamma^n}{n!}. \tag{32}$$

A small modification of the proof of Theorem 5 yields also the following generalization of (27), which is probably known but we have not been able to find a proof in the literature.

**Theorem 6.** (a) *Let $p>1$ and $a\in\mathcal{A}^m$, $m<-2n/p$. Then $W_a, T_a\in\mathcal{S}^p$.*
(b) *Let $p>1$ and $a\in\mathcal{A}^m$, $m\le -2n/p$. Then $W_a, T_a\in\mathcal{S}^{p,\infty}$.*
(c) *If $a\in\mathcal{A}^0$, then $W_a$ is bounded.*

*Proof.* We begin by showing how to reduce the assertions concerning $W_a$ to those about $T_a$. By Theorem 2 and (16), for $a\in\mathcal{A}^m$ with $m\le 0$,

$$T_a = W_a + W_{L_1 a} + W_{L_2 a} + \dots + W_{L_{2n-1}a} + W_{GLS^{m-2n}}, \tag{33}$$

where $L_k$ are the "operators" implicit in (20):

$$L_k f := \sum_{\substack{j,l\ge 0\\ j+2l=k}} \frac{\Delta^l}{l!(8\gamma)^l} a_j$$

for $a_j$ as in (19) (thus $L_k f\in\mathcal{A}^{m-k}$ if $f\in\mathcal{A}^m$); and $W_{GLS^{m-2n}}$ is a shorthand for "$W_g$ with some $g\in GLS^{m-2n}$". Replacing $a$ by $L_1 a$, we also have

$$T_{L_1 a} = W_{L_1 a} + W_{L_1(L_1 a)} + W_{L_2(L_1 a)} + \dots + W_{GLS^{m-2n}},$$

so

$$T_a = W_a + T_{L_1 a} + W_{(L_2 a - L_1 L_1 a)} + W_{(L_3 a - L_2 L_1 a)} + \cdots + W_{GLS^{m-2n}}.$$

Replacing $a$ in (33) by $L_2 a - L_1 L_1 a$ and continuing in this fashion, we get after $2n$ steps

$$T_a = W_a + T_{b_1} + T_{b_2} + \cdots + T_{b_{2n-1}} + W_{GLS^{m-2n}}$$

with $b_j \in \mathcal{A}^{m-j}$. Thus

$$W_a = T_b + W_{GLS^{m-2n}},$$

where $b \in \mathcal{A}^m$ has the same leading term in the expansion (19) as $a$.

For $m = 0$, $b$ is bounded, hence so is $T_b$, while $W_{GLS^{m-2n}} = W_{GLS^{-2n}}$ is Hilbert-Schmidt by (27). Thus $W_a$ is bounded, proving (c).

For $m < 0$, $W_{GLS^{m-2n}}$ is trace-class by (28), hence $W_a \in \mathcal{S}^p$ if and only if $T_b \in \mathcal{S}^p$ (for $p \geq 1$). Thus it is enough, as claimed, to prove (a), (b) only for $T_a$.

To achieve the latter — which again in fact hold even for any $a(z)$ which is $O(|z|^m)$, with $m$ as indicated, as $|z| \to +\infty$, i.e. no estimates on derivatives are needed — set, as in the proof of Theorem 5,

$$g(z) = (1 + |z|^2)^{-m/2} a(z),$$

so that $g \in L^\infty(\mathbf{C}^n)$ and

$$T_a = P_\gamma M_{(1+|z|^2)^{m/4}} M_g M_{(1+|z|^2)^{m/4}} P_\gamma.$$

Again, the fact that, for $B$ bounded, $A \in \mathcal{S}^{2p}$ and $C \in \mathcal{S}^{2p,\infty}$

$$A^* B A \in \mathcal{S}^p, \quad C^* B C \in \mathcal{S}^{p,\infty},$$

thus reduces the claims in (a), (b) to showing that the eigenvalues $c_j$ of $T_{(1+|z|^2)^{m/2}}$ satisfy

$$\sum_j c_j^p < \infty \quad \text{if } m < -\frac{2n}{p},$$

and

$$\sup_j (j+1)^{1/p} c_j < \infty \quad \text{if } m = -\frac{2n}{p},$$

respectively. However, a computation completely parallel to (29) exhibits that $T_{(1+|z|^2)^{m/2}}$ is diagonalized by the monomial basis $\{z^\alpha\}$, with eigenvalues

$$c_\alpha \leq \Big(\frac{\gamma}{|\alpha|}\Big)^{-m/2}.$$

Arguing as in (30), we then get

$$\sum_j c_j^p \leq \sum_k \binom{k+n-1}{k} \Big(\frac{\gamma}{k}\Big)^{-mp/2} \asymp \sum_k k^{n-1+\frac{mp}{2}}$$

which is finite for $m < -\frac{2n}{p}$; and

$$\sup_j (j+1)^{1/p} c_j \leq \sup_k \Big(\frac{k^n}{n!}\Big)^{1/p} \Big(\frac{\gamma}{k}\Big)^{-m/2} \asymp \sup_k k^{\frac{n}{p}+\frac{m}{2}},$$

which is finite if $m \leq -\frac{2n}{p}$. $\square$

## 5. Dixmier traces of Toeplitz and Weyl operators

The following result is a consequence of (32).

**Theorem 7.** *For $f \in \mathcal{A}^{-2n}$, $T_f$ and $W_f$ are measurable and*

$$\operatorname{tr}_\omega W_f = \operatorname{tr}_\omega T_f = \frac{\gamma^n}{n!} \int_{\mathbf{S}^{2n-1}} f_0(\zeta)\, d\sigma(\zeta),$$

*where $f_0$ is as in* (17) *and $d\sigma$ denotes the normalized surface measure on $\mathbf{S}^{2n-1}$.*

The proof below is adapted from [9].

*Proof.* We already know that $T_f, W_f \in \mathcal{S}^{\text{Dixm}}$ and $\operatorname{tr}_\omega T_f = \operatorname{tr}_\omega W_f$ for any $\omega$. Furthermore, for any $\chi \in C^\infty(\mathbf{R})$ that vanishes on the interval $(-\infty, \frac{1}{2})$ and equals 1 on $(1, +\infty)$, $f(z) - (1+|z|^2)^{-n}\chi(|z|)f_0(\frac{z}{|z|}) \in GLS^{-2n-2}$, so $W_f - W_{(1+|z|^2)^{-n}\chi f_0}$ is trace-class and has vanishing Dixmier trace; thus it is enough to prove the theorem only for functions $f$ of the form $f(z) = (1 + |z|^2)^{-n}\chi(|z|)f_0(\frac{z}{|z|})$, where $f_0 \in C^\infty(\mathbf{S}^{2n-1})$. Next, by (15), $f - (1 + |z|^2)^{-n} \# \chi(|z|)f_0(\frac{z}{|z|}) \in GLS^{-2n-2}$, so $W_f - W_{(1+|z|^2)^{-n}} W_{\chi(|z|)f_0(\frac{z}{|z|})}$ is trace-class and has vanishing Dixmier trace; and, finally, $W_{(1+|z|^2)^{-n}} - T_{(1+|z|^2)^{-n}}$ and $W_{f_0\chi} - T_{f_0\chi}$ are trace-class and in $\mathcal{S}^{n,\infty}$, respectively, by (27) and Theorem 6 (b). Altogether, we have

$$\operatorname{tr}_\omega T_f = \operatorname{tr}_\omega W_f = \operatorname{tr}_\omega (T_{(1+|z|^2)^{-n}} T_{\chi f_0}).$$

Consider now the linear functional on $C(\mathbf{S}^{2n-1})$

$$\Phi : C(\mathbf{S}^{2n-1}) \ni a \longmapsto \operatorname{tr}_\omega (T_{(1+|z|^2)^{-n}} T_{\chi a}) \in \mathbf{C}.$$

We know from Theorem 5 that $T_{(1+|z|^2)^{-n}} \in \mathcal{S}^{\text{Dixm}}$, while $T_{\chi a}$ is bounded with $\|T_{\chi a}\| \le \|\chi a\|_\infty = \|a\|_{C(\mathbf{S}^{2n-1})}$; thus by (25) $\Phi$ is bounded. By the Riesz representation theorem, there exists a measure $d\mu$ on $\mathbf{S}^{2n-1}$ such that

$$\Phi(a) = \int_{\mathbf{S}^{2n-1}} a\, d\mu.$$

For any unitary $n \times n$ matrix $U$, the "composition operator"

$$\tau_U : f(z) \mapsto f(Uz)$$

acts unitarily on $\mathcal{F}_\gamma$, and

$$\tau_U^* T_f \tau_U = T_{\tau_U f}$$

(in particular, $\tau_U^* T_{(1+|z|^2)^{-n}} \tau_U = T_{(1+|z|^2)^{-n}}$). From the cyclicity property (24) of the Dixmier trace we thus get $\Phi(\tau_U a) = \Phi(a)$ for all $a$, whence $d\mu(Uz) = d\mu(z)$. Thus $d\mu$ has to be invariant under all unitary transformations of $\mathbf{C}^n$; hence, it is a multiple of $d\sigma$, and

$$\operatorname{tr}_\omega (T_{(1+|z|^2)^{-n}} T_{\chi a}) = c \int_{\mathbf{S}^{2n-1}} a\, d\sigma$$

for some constant $c$ independent of $a$. Thus

$$\operatorname{tr}_\omega T_f = c \int_{\mathbf{S}^{2n-1}} f_0(\zeta)\, d\sigma(\zeta).$$

Taking in particular $f(z) = (1 + |z|^2)^{-n}$ and using (32) gives $c = \gamma^n / n!$. □

## 6. Dixmier traces of Hankel operators

From the formula

$$H_f^* H_g = T_{\overline{f}g} - T_{\overline{f}} T_g$$

and (16) we see that

$$\begin{aligned} H_f^* H_g &= W_{\mathcal{E}(\overline{f}g)} - W_{\mathcal{E}\overline{f}} W_{\mathcal{E}g} \\ &= W_{\mathcal{E}(\overline{f}g) - \mathcal{E}\overline{f} \# \mathcal{E}g}. \end{aligned}$$

Assuming that $f \in \mathcal{A}^m$, $g \in \mathcal{A}_k$, $m, k \le 0$, let us temporarily denote by $O(-4)$ a general function in $\mathcal{A}^{m+k-4}$. By (20) and (15), we then have

$$\begin{aligned} \mathcal{E}(\overline{f}g) &= \overline{f}g + \frac{\Delta(\overline{f}g)}{8\gamma} + O(-4), \\ \mathcal{E}\overline{f} \# \mathcal{E}g &= \Big(\overline{f} + \frac{\Delta\overline{f}}{8\gamma} + O(-4)\Big) \# \Big(g + \frac{\Delta g}{8\gamma} + O(-4)\Big) \\ &= \overline{f} \# g + \frac{\overline{f} \# \Delta g + \Delta\overline{f} \# g}{8\gamma} + O(-4) = \overline{f} \# g + \frac{\overline{f}\Delta g + g\Delta\overline{f}}{8\gamma} + O(-4), \end{aligned}$$

and

$$\begin{aligned} \mathcal{E}(\overline{f}g) - \mathcal{E}\overline{f} \# \mathcal{E}g &= \overline{f}g - \overline{f} \# g + \frac{\Delta(\overline{f}g) - \overline{f}\Delta g - g\Delta\overline{f}}{8\gamma} + O(-4) \\ &= \frac{1}{2\gamma} \sum_{|\alpha|+|\beta|=1} \frac{(-1)^{|\beta|}}{\alpha!\beta!} \partial^\alpha \overline{\partial}^\beta f \cdot \partial^\beta \overline{\partial}^\alpha g + \frac{\Delta(\overline{f}g) - \overline{f}\Delta g - g\Delta\overline{f}}{8\gamma} + O(-4) \\ &= \frac{1}{2\gamma} \sum_{j=1}^n (\partial_j \overline{f} \overline{\partial}_j g - \overline{\partial}_j \overline{f} \partial_j g) + \frac{1}{2\gamma} \sum_{j=1}^n (\partial_j \overline{f} \overline{\partial}_j g + \overline{\partial}_j \overline{f} \partial_j g) + O(-4) \\ &= \frac{1}{\gamma} \sum_{j=1}^n \partial_j \overline{f} \overline{\partial}_j g + O(-4), \end{aligned}$$

where we have used the Leibniz rule

$$\frac{\Delta}{4}(\overline{f}g) = \sum_{j=1}^n \partial_j \overline{\partial}_j (\overline{f}g) = \sum_{j=1}^n (\partial_j \overline{\partial}_j f \cdot g + \partial_j f \cdot \overline{\partial}_j g + \overline{\partial}_j f \cdot \partial_j g + f \cdot \partial_j \overline{\partial}_j g).$$

We have thus proved the following.

**Proposition 8.** *For $f \in \mathcal{A}^m$, $g \in \mathcal{A}^k$, $m, k \le 0$, we have*

$$\mathcal{E}(\overline{f}g) - \mathcal{E}\overline{f} \# \mathcal{E}g = \frac{1}{\gamma} \sum_{j=1}^n \partial_j \overline{f} \cdot \overline{\partial}_j g + h, \qquad h \in \mathcal{A}^{m+k-4}.$$

**Corollary 9.** *For $f,g \in \mathcal{A}^0$ and $\zeta \in \mathbf{S}^{2n-1}$, denote*

$$Q(f,g)(\zeta) := \lim_{r\to+\infty} r^2 \sum_{j=1}^n \partial_j \overline{f}(r\zeta)\cdot \overline{\partial}_j g(r\zeta) \tag{34}$$

(*the limit exists in view of* (17)).

*Then $(H_f^* H_g)^n$ belongs to $\mathcal{S}^{\rm Dixm}$, is measurable, and*

$$\operatorname{tr}_\omega (H_f^* H_g)^n = \frac{1}{n!}\int_{\mathbf{S}^{2n-1}} Q(f,g)^n \, d\sigma. \tag{35}$$

*Proof.* By (13) and the preceding proposition, $(H_f^* H_g)^n = W_a$ where $a = (\mathcal{E}(\overline{f}g) - \mathcal{E}\overline{f} \# \mathcal{E}g)^{\# n} \in \mathcal{A}^{-2n}$ satisfies

$$a = \Big(\frac{1}{\gamma}\sum_{j=1}^n \partial_j \overline{f}\, \overline{\partial}_j g\Big)^n + h, \qquad h \in \mathcal{A}^{-2n-2}.$$

The membership of $(H_f^* H_g)^n$ in $\mathcal{S}^{\rm Dixm}$ and the formula for its Dixmier trace thus follow by Theorem 7. Since the right-hand side of (35) does not depend on $\omega$, $(H_f^* H_g)^n$ must be measurable. □

Similarly, from the relation

$$[T_f, T_g] = (T_{fg} - T_g T_f) - (T_{fg} - T_f T_g) = H_{\overline{g}}^* H_f - H_{\overline{f}}^* H_g$$

we immediately get:

**Corollary 10.** *For $f_1, g_1, \dots, f_n, g_n \in \mathcal{A}^0$, the operator $T = [T_{f_1}, T_{g_1}] \dots [T_{f_n}, T_{g_n}]$ belongs to $\mathcal{S}^{\rm Dixm}$, is measurable, and*

$$\operatorname{tr}_\omega T = \frac{1}{n!}\int_{\mathbf{S}^{2n-1}} \prod_{j=1}^n (Q(\overline{g}, f) - Q(\overline{f}, g)) \, d\sigma,$$

*with $Q(f,g)$ as in* (34).

The same argument yields also the following more general result.

**Corollary 11.** *Let $f_1 \in \mathcal{A}^{-m_1}$, $g_1 \in \mathcal{A}^{-k_1}$, …, $f_l \in \mathcal{A}^{-m_l}$, $g_l \in \mathcal{A}^{-k_l}$, $m_1, \dots, m_l$, $k_1, \dots, k_l \geq 0$, be such that $m_1 + \dots + m_l + k_1 + \dots + k_l + 2l = 2n$, and let $h_1, \dots, h_s \in \mathcal{A}^0$. Denote, for $j = 1, \dots, s$ and $\zeta \in \mathbf{S}^{2n-1}$,*

$$P(h_j)(\zeta) := \lim_{r\to+\infty} h_j(r\zeta),$$

*and for $j = 1, \dots, l$ and $\zeta \in \mathbf{S}^{2n-1}$,*

$$Q_{m_j k_j}(f_j, g_j)(\zeta) := \lim_{r\to+\infty} r^{m_j + k_j + 2} \sum_{j=1}^n \partial_j \overline{f}(r\zeta)\cdot \overline{\partial}_j g(r\zeta)$$

(*the existence of the limits is again guaranteed by* (17)). *Then the product of* $(H^*_{f_1}H_{g_1}),\dots,(H^*_{f_l}H_{g_l}),T_{h_1},\dots,T_{h_s}$ (*in any order*) *belongs to* $\mathcal{S}^{\rm Dixm}$, *is measurable, and its Dixmier trace equals*

$$\frac{\gamma^{n-l}}{n!}\int_{\mathbf{S}^{2n-1}} Q_{m_1k_1}(f_1,g_1)\dots Q_{m_lk_l}(f_l,g_l)P(h_1)\dots P(h_s)\,d\sigma. \tag{36}$$

Clearly, for $f$ with the expansion (17), $P(f)=f_0$, the leading order term. We conclude this section by deriving a formula expressing $Q(f,g)$ in terms of tangential derivatives on the unit sphere of the leading order terms $f_0,g_0$ of $f,g$, thus completing the proof of Theorem 1.

Let

$$R=\sum_{j=1}^n z_j\partial_j, \qquad \overline{R}=\sum_{j=1}^n \overline{z}_j\overline{\partial}_j,$$

denote the holomorphic and anti-holomorphic components, respectively, of the radial derivative. The vector fields

$$\partial^b_j := \partial_j-\frac{\overline{z}_j}{|z|^2}R, \quad \overline{\partial}^b_j := \overline{\partial}_j-\frac{z_j}{|z|^2}\overline{R}, \tag{37}$$

defined on $\mathbf{C}^n\setminus\{0\}$, are tangential to the spheres $|z|=\text{const}$. In view of the linear dependencies

$$\sum_{j=1}^n z_j\partial^b_j=\sum_{j=1}^n \overline{z}_j\overline{\partial}^b_j=0, \tag{38}$$

they span at each $z\neq0$ a subspace of complex dimension $2n-2$. Together with the Reeb vector field

$$E=\frac{\overline{R}-R}{2|z|}, \tag{39}$$

which is also tangent to $|z|=\text{const.}$, they span the complexified tangent space to the spheres at each $z\neq0$. The orthogonal complement is given by the normal differentiation operator

$$N=\frac{R+\overline{R}}{|z|}=\frac{\partial}{\partial r} \tag{40}$$

for $z=r\zeta$, $r>0$, $\zeta\in\mathbf{S}^{2n-1}$. Note that all the vector fields (37)–(40) have been normalized to lower the homogeneity degree of a function by one, just as the ordinary derivatives $\partial_j,\overline{\partial}_j$ do.

**Proposition 12.** *For* $f\in\mathcal{A}^{-m}$, $g\in\mathcal{A}^{-k}$, $m,k\geq0$, *with asymptotic expansions*

$$f(z)\sim\sum_{j=m}^\infty |z|^{-j}f_j\Big(\frac{z}{|z|}\Big), \quad g(z)\sim\sum_{j=k}^\infty |z|^{-j}g_j\Big(\frac{z}{|z|}\Big), \quad |z|\geq1,$$

*we have*

$$Q_{mk}(f,g)=\Big(\frac{m}{2}+E\Big)f_m\cdot\Big(\frac{k}{2}-E\Big)g_k+\sum_{j=1}^n\partial_j^b\overline{f}_m\cdot\overline{\partial}_j^b g_k.$$

*Proof.* Looking at the homogeneity degree it is immediate that $Q_{mk}(f,g)$ depends only on the top order terms of $f,g$, so without loss of generality we may assume that in fact $f(z)=|z|^{-m}f_m(\frac{z}{|z|})$ and $g(z)=|z|^{-k}g_k(\frac{z}{|z|})$ for $|z|>\frac12$. Denote momentarily, for brevity, $F(z)=f_m(\frac{z}{|z|})$. From the equality

$$F(z)=F(rz),\qquad \forall r>0,$$

we get by the chain rule

$$\partial_j^b F(z)=r\partial_j^b F(rz),$$

hence for $\zeta\in\mathbf{S}^{2n-1}$,

$$\partial_j^b F(r\zeta)=\frac{\partial_j^b F(\zeta)}{r}=\frac{\partial_j^b f_m(\zeta)}{r},$$

since, $\partial_j^b$ being tangential to $\mathbf{S}^{2n-1}$, the derivative $\partial_j^b F$ depends only on the restriction $f_m$ of $F$ to $\mathbf{S}^{2n-1}$. As $\partial_j^b$ annihilates $|z|$, we thus get

$$\partial_j^b f(z)=\frac{\partial_j^b F(z)}{|z|^m}=\frac{\partial_j^b f_m(\frac{z}{|z|})}{|z|^{m+1}},$$

that is,

$$\partial_j^b f(z)=\frac{\partial_j^b f_m(\zeta)}{r^{m+1}} \tag{41}$$

for $z=r\zeta$, $r>\frac12$, $\zeta\in\mathbf{S}^{2n-1}$. Similarly

$$\overline{\partial}_j^b f(r\zeta)=\frac{\overline{\partial}_j^b f_m(\zeta)}{r^{m+1}},\quad Ef(r\zeta)=\frac{Ef_m(\zeta)}{r^{m+1}}, \tag{42}$$

while

$$Nf(r\zeta)=-\frac{mf_m(\zeta)}{r^{m+1}} \tag{43}$$

since $N$ annihilates $F$.

From the definitions (37)–(40), we have

$$\partial_j=\partial_j^b+\frac{\overline{z}_j}{|z|}\Big(\frac{N}{2}-E\Big),\qquad \overline{\partial}_j=\overline{\partial}_j^b+\frac{z_j}{|z|}\Big(\frac{N}{2}+E\Big).$$

Thus

$$r^{m+k+2}\sum_{j=1}^n\partial_j\overline{f}\cdot\overline{\partial}_j g=\sum_{j=1}^n\Big(\partial_j^b\overline{f}_m-\overline{\zeta}_j\Big(\frac{m}{2}\overline{f}_m+E\overline{f}_m\Big)\Big)\Big(\overline{\partial}_j^b g_k-\zeta_j\Big(\frac{k}{2}g_k-Eg_k\Big)\Big).$$

Multiplying out and using the relations (38) and $|\zeta|^2 = 1$ gives

$$r^{m+k+2}\sum_{j=1}^n \partial_j\overline{f}\cdot\overline{\partial}_j g = \sum_{j=1}^n \partial_j^b\overline{f}_m\overline{\partial}_j^b g_k + \Big(\frac{m}{2}\overline{f}_m + E\overline{f}_m\Big)\Big(\frac{k}{2}g_k - Eg_k\Big),$$

completing the proof of the claim. □

Taking in particular $m = k = 0$, we get

$$Q_{00}(f,g) \equiv Q(f,g) = \{\overline{f}_0, g_0\}_F,$$

where

$$\{\phi,\psi\}_F := \sum_{j=1}^n \partial_j^b\phi\cdot\overline{\partial}_j^b\psi - E\phi\cdot E\psi, \tag{44}$$

thereby proving Theorem 1.

## 7. Concluding remarks

**7.1 Spherical Laplacian.** A short computation gives

$$\frac{\Delta}{4} = \sum_{j=1}^n \partial_j\overline{\partial}_j = \sum_{j=1}^n \partial_j^b\overline{\partial}_j^b - E^2 + \frac{N^2}{4},$$

meaning that

$$L := \sum_{j=1}^n \partial_j\overline{\partial}_j = \sum_{j=1}^n \partial_j^b\overline{\partial}_j^b - E^2$$

is the Laplace operator on the sphere $\mathbf{S}^{2n-1}$.

Note that, in some sense, (44) is obtained from $L$ upon "applying holomorphic derivatives to $\phi$, and anti-holomorphic ones to $\psi$".

**7.2 Separation of variables.** Toeplitz operators $T_f$ on $\mathcal{F}_\gamma$ make sense also for suitable unbounded symbols $f$, e.g. for $f$ a polynomial in $z, \overline{z}$, as densely defined operators. For $f$ holomorphic, since a product of two holomorphic functions is again holomorphic, it follows that $T_f$ is just the operator of "multiplication by $f$" on $\mathcal{F}_\gamma$, from which one gets

$$T_{\overline{f}}T_g = T_{\overline{f}g}, \quad T_gT_f = T_{fg}, \qquad \forall g \in L^\infty(\mathbf{C}^n).$$

In particular, $Q(f,g)$ must vanish if $f$ or $g$ is holomorphic. This explains why there are only holomorphic derivatives of $\overline{f}$ and anti-holomorphic derivatives of $g$ in (34).

**7.3 Still other generalizations.** Using the same ideas, it is possible to prove many variants of Theorem 7 and Corollary 11; for instance, for any $f_1 \in \mathcal{A}^{-m_1}$, ..., $f_s \in \mathcal{A}^{-m_s}$, $m_1, \dots, m_s \ge 0$, with $m_1 + \cdots + m_s = 2n$, the product $T_{f_1}\dots T_{f_s}$ belongs to $\mathcal{S}^{\mathrm{Dixm}}$, is measurable and

$$\operatorname{tr}_\omega T_{f_1}\dots T_{f_s} = \frac{\gamma^n}{n!}\int_{\mathbf{S}^{2n-1}} \lim_{r\to+\infty} r^{2n} f_1(r\zeta)\dots f_s(r\zeta)\,d\sigma(\zeta), \tag{45}$$

and similarly for $W_{f_1}\dots W_{f_s}$.

**7.4 More general symbols.** It would be interesting to know whether (45) remains in force for more general functions $f_1, \dots, f_s$ (e.g. not necessarily belonging to the $\mathcal{A}^m$ classes), possibly with the "integral of limit" replaced by "limit of integral", provided the latter exists.

**7.5 Size does not suffice.** In the proofs of Theorems 5 and 6, we in fact showed that $T_f \in \mathcal{S}^p$, $1 \leq p < \infty$, not only for $f \in \mathcal{A}^m$ with $m < -\frac{2n}{p}$, but even for any $f$ which is $O(|z|^{-m})$, as $|z| \to +\infty$, with such $m$; and similarly for $\mathcal{S}^{p,\infty}$, while obviously $T_f$ is bounded if $f$ is bounded. This is no longer true for $W_f$; indeed, there exist bounded smooth $f$ for which $W_f$ fails to be a bounded operator [23, Chapter 12]. However, the hypothesis $f \in \mathcal{A}^m$ can be relaxed to $f \in GLS^m$ (for Theorems 2, 5, 6, and Proposition 8) or at least to $f \in \mathcal{A}^m + GLS^{m-\delta}$, $\delta > 0$ (for Theorem 7, Corollaries 9, 10 and Proposition 11).

**7.6 Positive orders.** It is not difficult to extend some, if not all, of our results also to $f \in \mathcal{A}^m$ with positive $m$ (the Toeplitz operators $T_f$ being then only densely defined in general). The only point needing some adjustments is the proof of the convolution estimate in Lemma 3 for $m > 0$, which is elementary.

**7.7 Comparison with the ball.** Our results should be compared to those for Toeplitz/Hankel operators on the unit ball $\mathbf{B}^n$ of $\mathbf{C}^n$ in [11] (cf. (4) and (7)), obtained via the "pseudo-Toeplitz operators" of Howe [18], as recalled in the Introduction: the "boundary Poisson bracket" $\{f, g\}_b$ there, which is the counterpart of our $\{g, f\}_F - \{f, g\}_F$, is actually given by the same formula as the latter; however, the "half-bracket" $\{g, \overline{f}\}_L$ (corresponding to the symbol of the product $H^*_{\overline{f}} H_g$, rather than of the commutator $[T_f, T_g] = H^*_{\overline{g}} H_f - H^*_{\overline{f}} H_g$), playing the role of our $\{f, g\}_F$, differs from (44) by the absence of the term involving the Reeb vector field $E$ (cf. Section 4 in [11] and p. 619 in [12]).

(H.B. and H.Y.) LATP, U.M.R. C.N.R.S. 6632, CMI, Université de Provence, 39 Rue F-Joliot-Curie, 13453 Marseille Cedex 13, France
*E-mail address*: bommier@gyptis.univ-mrs.fr, youssfi@gyptis.univ-mrs.fr

(M.E.) Mathematics Institute, Žitná 25, 11567 Prague 1, Czech Republic and Mathematics Institute, Silesian University at Opava, Na Rybníčku 1, 74601 Opava, Czech Republic
*E-mail address*: englis@math.cas.cz